\documentclass[12pt]{article}
\usepackage{amsfonts}
\usepackage{amsmath}
\usepackage{mathrsfs}
\usepackage{color}
\usepackage{tikz}
\usepackage{authblk}
\usepackage{mathrsfs,amscd,amssymb,amsthm,amsmath,bm,graphicx,psfrag,subfigure,url,verbatim}

\allowdisplaybreaks
\usepackage{indentfirst}
\newtheorem{thm}{Theorem}

\newtheorem{claim}{Claim}

\newtheorem{lem}{Lemma}
\newtheorem{cor}{Corollary}

\newenvironment{wst}
{\setlength{\leftmargini}{1.5\parindent}
 \begin{itemize}
 \setlength{\itemsep}{-1.1mm}}
{\end{itemize}}

\begin{document}
	
\title{\bf Pairwise Intersection Profiles of Intersecting Families}
\author[a]{Sumin Huang\thanks{Email: sumin@nuist.edu.cn}}
\affil[a]{School of Mathematics and Statistics, Nanjing University of Information Science $\&$ Technology, Nanjing 210044, P. R. China}
\date{}
\maketitle

\begin{abstract}
Let $\mathcal F\subseteq\binom{[n]}{k}$ be a $k$-uniform family, and let $P_m(\mathcal F)$ denote the number of ordered pairs $(A,B)\in\mathcal F\times\mathcal F$ with $|A\cap B|=m$. We study the pairwise intersection profile $(P_0(\mathcal F),P_1(\mathcal F),\ldots,P_k(\mathcal F))$ of intersecting families. We prove that for every $1\leq m\leq k$, if $n\geq\max\{2k,3k-m-1\}$, then 
$$
P_m(\mathcal F)\leq P_m(\mathcal S_1^k),
$$
where $\mathcal S_1^k$ is a full star. Moreover, if $n\geq \max\{2k+1,3k-m\}$, then equality holds if and only if $\mathcal F$ is a star up to permutations. As a corollary, for any function $f:\{0,1,\ldots,k\}\rightarrow \mathbb{R}_{\geq 0}$ and $n\geq 3k-2$, we show that $\Phi_f(\mathcal F)=\sum_{(A,B)\in\mathcal F\times\mathcal F}f(|A\cap B|)$ is maximized by a star. We also establish a stability result showing that, for $n\geq \max\{2k+1,3k-m\}$, if $P_m(\mathcal F)$ is sufficiently close to its maximum, then $\mathcal F$ has an element of degree close to $\binom{n-1}{k-1}$.
\end{abstract}

\vspace{2mm} \noindent{\bf Keywords}: pairwise intersection profiles, intersecting family, stability, Katona's circle method
\vspace{2mm}

\setcounter{section}{0}

\section{Introduction}
  For integers $n\geq k\geq 0$, let $[n]=\{1,\ldots,n\}$ and let $\binom{[n]}{k}$ denote the family of all $k$-subsets of $[n]$. A family $\mathcal{F}\subseteq \binom{[n]}{k}$ is called {\it intersecting} if $A\cap B\neq \emptyset$ for any $A,B\in\mathcal{F}$. In particular, if all members of $\mathcal{F}$ contain a fixed element, then $\mathcal{F}$ is called a {\it star}, and is trivially intersecting. For convenience, let $\mathcal{S}^k_1:=\{S\in \binom{[n]}{k}:1\in S\}$ denote the star centered at $1$. 
  
   The celebrated Erd\H{o}s--Ko--Rado theorem \cite{Erdos1961} states that $\mathcal{S}^k_1$ has the largest size among all intersecting families of $k$-sets provided $n\geq 2k$:
    \begin{thm}\cite{Erdos1961}\label{thm-ekr}
    	Let $\mathcal F\subseteq\binom{[n]}{k}$ be an intersecting family. If $n\ge2k$, then $|\mathcal{F}|\leq \binom{n-1}{k-1}.$ 
    	Moreover, if $n>2k$, equality holds if and only if $\mathcal{F}=\mathcal{S}^k_1$ up to permutations.
    \end{thm}
   The Erd\H{o}s--Ko--Rado theorem has led to many important developments in extremal set theory, including the Hilton--Milner theorem \cite{Hilton1967}, the $r$-wise intersection theorem \cite{Frankl1976}, and the Complete Intersection Theorem \cite{Ahlswede1997}. More generally, determining the maximum size of a set family under various intersection restrictions has been one of the central problems in extremal set theory.
    
   The size of a family, however, does not describe how its members intersect with each other. To record the pairwise intersection structure more precisely, for a family $\mathcal F\subseteq\binom{[n]}{k}$ and $0\leq m\leq k$, define
   $$
   P_m(\mathcal F):=\left|\left\{(A,B)\in\mathcal{F}\times\mathcal{F}:|A\cap B|=m\right\}\right|.
   $$
   Thus, $P_m(\mathcal F)$ denotes the number of ordered pairs of sets of $\mathcal F$ whose intersection has size exactly $m$. We call
   $$
   \textbf{P}(\mathcal{F})=(P_0(\mathcal F),P_1(\mathcal F),\ldots,P_k(\mathcal F))
   $$
   {\it the pairwise intersection profile} of $\mathcal{F}$. Clearly, $P_k(\mathcal{F})=|\mathcal{F}|$. Also, $\mathcal{F}$ is intersecting if and only if $P_0(\mathcal{F})=0$. Consequently, the Erd\H{o}s--Ko--Rado theorem can be viewed as determining the maximum value of $P_k(\mathcal{F})$ provided $P_0(\mathcal{F})=0$. 
    
   Several classical directions in extremal set theory can be naturally interpreted in terms of the pairwise intersection profile. Ray--Chaudhuri and Wilson \cite{Ray-Chaudhuri1975} determined upper bounds on the size of $L$-intersecting families. Equivalently, they determined upper bounds on $P_k(\mathcal{F})$ when $P_m=0$ for every $m\notin [k-1]\backslash L$. Similarly, the forbidden-intersection problem prescribes a single vanishing coordinate $P_t(\mathcal{F})=0$ and asks for the maximum possible value of $P_k(\mathcal{F})$. Important results in this direction were obtained by Frankl and F\"uredi \cite{Frankl1985}, while the case $t=1$ was further studied by Keevash, Mubayi and Wilson \cite{Keevash2006}. As a supersaturation problem, Das, Gan and Sudakov \cite{Das2016} studied the minimum possible value of $P_0(\mathcal F)$, namely the number of disjoint pairs, among families with given $P_k(\mathcal{F})$. Recently, Wu and Zhang \cite{Wu2026} determined the maximum value of $P_{k-1}(\mathcal{F})$ for $t$-intersecting families $\mathcal{F}$, which in terms of the intersection profile are precisely those satisfying $P_0(\mathcal{F})=P_1(\mathcal{F})=\cdots=P_{t-1}(\mathcal{F})=0$.

   The preceding problems suggest that many classical questions in extremal set theory can be viewed as fixing or constraining certain coordinates of the pairwise intersection profile and asking for the extremal values of other coordinates. In particular, under the basic constraint $P_0(\mathcal{F})=0$, or equivalently $\mathcal{F}$ is intersecting, it is natural to ask for the maximum possible value of each coordinate $P_m(\mathcal F)$.
   
   In this paper, our main result shows that, for any $m\in [k]$, when $n\geq\max\{2k,3k-m-1\}$, $P_m(\mathcal{F})$ is maximized by a star among all intersecting families $\mathcal F\subseteq\binom{[n]}{k}$.
   
   \begin{thm}\label{thm-pm}
   	Let
   	$\mathcal F\subseteq\binom{[n]}{k}$ be an intersecting family. Then,
   	for every $m\in [k]$, if $n\geq\max\{2k,3k-m-1\}$, then
   	$$
   	P_m(\mathcal F)\leq P_m(\mathcal S_1^k).
   	$$
   	Moreover, if $n\geq \max\{2k+1,3k-m\}$, then equality holds if and only if $\mathcal{F}=\mathcal{S}_1^k$ up to permutations.
   \end{thm} 
    
   When $m=k$, Theorem~\ref{thm-pm} reduces to the Erd\H{o}s--Ko--Rado theorem.
   The equality case in Theorem~\ref{thm-pm} further shows that attaining the extremal value at any one level determines the entire pairwise intersection profile.
    
   \begin{cor}
   	Let $n\geq 2k+1$ and
   	$\mathcal F\subseteq\binom{[n]}{k}$ be an intersecting family. If there
   	exists some $m\in[k]$ such that $n\geq 3k-m$ and
   	$$
   	P_m(\mathcal F)=P_m(\mathcal S_1^k),
   	$$
   	then
   	$$
   	P_i(\mathcal F)=P_i(\mathcal S_1^k)
   	$$
   	for every $0\leq i\leq k$. In particular, $\mathcal F$ is a star up to
   	permutations.
   \end{cor}

   Theorem~\ref{thm-pm} determines the extremal behavior of each level of the pairwise intersection profile. A natural consequence is that the star also maximizes any weighted sum of these parameters. More precisely, for a function $f:\{0,1,\ldots,k\}\rightarrow\mathbb R$, define
   $$
   \Phi_f(\mathcal F):=\sum_{m=0}^{k}f(m)P_m(\mathcal F).
   $$
   Equivalently,
   $$
   \Phi_f(\mathcal F)=\sum_{(A,B)\in\mathcal F\times\mathcal F}f(|A\cap B|).
   $$
   The parameter $\Phi_f(\mathcal F)$ provides a flexible way to measure the overall contribution of intersections in a family by assigning different weights to different intersection sizes. By Theorem~\ref{thm-pm}, the following corollary holds immediately.
    
   \begin{cor}\label{cor-phi}
   	Let
   	$\mathcal F\subseteq\binom{[n]}{k}$ be an intersecting family. Let $f:\{0,1,\ldots,k\}\rightarrow \mathbb{R}_{\geq 0}$ be a nonnegative function. If $n\geq 3k-2$, then
   	$$
   	\Phi_f(\mathcal F)\leq\Phi_f(\mathcal S_1^k).
   	$$
   	Moreover, if $f(m)>0$ for some $m\in[k]$ such that $n\geq \max\{2k+1,3k-m\}$, then equality
   	holds if and only if
   	$\mathcal F=\mathcal S_1^k$ up to permutations.
   \end{cor} 

   Various choices of the weight function $f$ lead to intersection parameters that have been studied previously \cite{Huang2026,Ge2026,Wu2026}; see Section~5 for details.
   
   We also prove the following stability theorem. Let $d_{\mathcal F}(x)$ be the degree of $x$ in $\mathcal{F}$. 

   \begin{thm}\label{thm-stability}
   	Let $m\in [k]$ and $n\geq \max\{2k+1,3k-m\}$.
   	For every $\varepsilon>0$, there exists  $0<\delta<\min\{\frac{1}{m(n^4+1)},\frac{\varepsilon}{m(n^3+1)}\}$ such that if
   	$\mathcal F\subseteq\binom{[n]}{k}$ is an intersecting family satisfying
   	$P_m(\mathcal F)\geq (1-\delta)P_m(\mathcal S_1^k),$
   	then there exists $x\in[n]$ such that
   	$$
   	d_{\mathcal F}(x)\geq
   	(1-\varepsilon)\binom{n-1}{k-1}.
   	$$
   \end{thm}   
   
   The rest of the paper is organized as follows. In Section~2, we introduce cyclic permutations and several preliminary lemmas. In Section~3, we prove the upper bound in Theorem~\ref{thm-pm} using Katona's circle method \cite{Katona1972}. In Section~4, we prove Theorem~\ref{thm-stability} via an expansion argument on an adjacent-transposition Cayley graph, following the stability approach of Keevash \cite{Keevash2008} and Kamat \cite{Kamat2011}. We then use Theorem~\ref{thm-stability} to complete the proof of the equality case in Theorem~\ref{thm-pm}. Finally, in Section~5, we give some concluding remarks, including a discussion of the range of $n$ and connections between special choices of the weight function $f$ and previously studied intersection parameters.

   At the end of this section, it is worth mentioning that Talbot \cite{Talbot2004,Talbot2005} studied a related but different aspect of the intersection structure. Instead of counting the number of pairs with a prescribed intersection size, Talbot considered the number of distinct intersections, that is the size of $\{A\cap B:A,B\in \mathcal{F},|A\cap B|=m\}$. This quantity counts the possible $m$-intersections, whereas $P_m(\mathcal F)$ counts the number of pairs producing such intersections.

\section{Preliminaries}
In this section, we introduce some necessary definitions and lemmas. 

Consider a permutation $\sigma\in S_n$ as a sequence $(\sigma(1),\sigma(2),\ldots,\sigma(n))$. We say two permutations $\sigma_1$ and $\sigma_2$ are equivalent if there exists $i\in [n]$ such that $\sigma_1(x)=\sigma_2(x+i)$ for all $x\in [n]$, where the indices are taken modulo $n$. Let $C_n$ be a system of representatives of the equivalence
classes. We call the elements of $C_n$ {\it cyclic permutations} of $[n]$.
It is clear that each $\pi\in C_n$ represents a unique cyclic order $(\pi(1),\pi(2),\ldots,\pi(n))$ of $[n]$, and $|C_n|=(n-1)!$. For a cyclic order $\pi\in C_n$, an $r$-interval of $\pi$ is a set of the form $\{\pi(i),\pi(i+1),\ldots,\pi(i+r-1)\}$ for some $i\in[n]$. 

Using cyclic order, Katona proved the following lemma, which gives a nice simple proof of Theorem~\ref{thm-ekr}.
\begin{lem}\cite{Katona1972}\label{lem-interval}
	Let $\pi$ be a cyclic permutation of $[n]$ and let $G_1,\ldots,G_r$ be $k$-intervals of $\pi$ that form an intersecting family. Then $r\leq k$ provided $n\geq 2k$. Furthermore, if $n>2k$ and $r=k$, then $G_1,\ldots,G_r$ are all the $k$-intervals that contain a fixed element $x$.
\end{lem}

Note that $C_n$ is isomorphic to a symmetric group $\mathrm{Sym}([n]\backslash \{c\})$ for each $c\in [n]$. The following lemma was shown by Keevash \cite{Keevash2008}, using a result of Bacher \cite{Bacher1994}. We consider the following $(n-2)$-regular Cayley graph $G_c$. 
\begin{lem}\cite{Keevash2008}\label{lem-cay}
	Let $c\in [n]$ and $G_c$ be the Cayley graph on $\mathrm{Sym}([n]\backslash \{c\})$ generated by the set of adjacent transpositions 
	$$
	A=\{(1~2),(2~3),\ldots,(c-2~c-1),(c+1~c+2),\ldots,(n-1~n),(n~1)\}.
	$$
	That is $V(G_c)=\mathrm{Sym}([n]\backslash \{c\})$ and two permutations $\pi$ and $\mu$ are adjacent if and only if $\pi=\mu\sigma$ for some $\sigma\in A$.
	Then for any $H\subseteq V(G_c)$ with $|H|\leq \frac{|V(G_c)|}{2}$, we have
	$$
	|N_{G_c}(H)|\geq \frac{|H|}{n^3},
	$$
	where $N_{G_c}(H)$ is the set of vertices in $V(G_c)\backslash H$ which are adjacent to some vertex in $H$.
\end{lem}

\section{Proof of the Upper Bound in Theorem~\ref{thm-pm}}
In this section, we prove the upper bound in Theorem~\ref{thm-pm}. The characterization of the equality case will be deduced from the stability result in the next section.

If $m=k$, then $P_{m}(\mathcal{F})=|\mathcal{F}|$ and Theorem~\ref{thm-ekr} implies the desired result. Now we assume that $m<k$.

For a pair of sets $(A,B)\in \mathcal{F}\times\mathcal{F}$ and a fixed cyclic permutation $\pi$, we say $(A,B)$ is a \textit{representable pair} of $\pi$ if the following conditions hold:
\begin{wst}
	\item $A$ and $B$ are distinct intervals of $\pi$;
	\item the intersection $A\cap B$ and $A$ have the same right endpoint;
	\item the intersection $A\cap B$ and $B$ have the same left endpoint.
\end{wst}
Let $\mathcal{R}_{\pi}$ denote the set of all representable pairs of $\pi$ and $\mathcal{F}_{\pi}$ denote the set of intervals of $\pi$ that belong to $\mathcal{F}$.

Let $\mathcal{P}_m:=\{(A,B)\in \mathcal{F}\times \mathcal{F}:|A\cap B|=m\}$. 
We double count the size of the following set 
$$X_m:=\{(\pi,A,B):\pi\in C_n, (A,B)\in \mathcal{P}_m\cap \mathcal{R}_{\pi}\}.$$

For a fixed pair $(A,B)\in\mathcal{P}_m$, there exist exactly $m!(k-m)!^2(n-2k+m)!$ cyclic permutations $\pi$ such that $(A,B)$ is a representable pair of $\pi$. Thus 
\begin{align}\label{eq3}
	|X_m|=\sum_{(A,B)\in\mathcal{P}_m}m!(k-m)!^2(n-2k+m)!=m!(k-m)!^2(n-2k+m)!P_m(\mathcal{F}).
\end{align}

On the other hand, for a fixed $\pi\in C_n$, let 
$$
\mathcal{I}_{m,\pi}:= \{A\cap B: (A,B)\in \mathcal{P}_m\cap \mathcal{R}_{\pi}\}.
$$

\begin{claim}\label{cla-m}
$\mathcal{I}_{m,\pi}$ is an intersecting family consisting of $m$-intervals.
\end{claim}
\begin{proof}[Proof of Claim~\ref{cla-m}]
Suppose that there exist two disjoint intervals $A\cap B$ and $A'\cap B'$ in $\mathcal{I}_{m,\pi}$. Without loss of generality, assume that $A\cap B=\{\pi(1),\ldots,\pi(m)\}$ and $A'\cap B'=\{\pi(i),\ldots,\pi(i+m-1)\}$ for some $m+1\leq i\leq n-m+1$. Then 
\begin{align*}
	A=&\{\pi(n-k+m+1),\pi(n-k+m+2),\ldots,\pi(n),\pi(1),\ldots,\pi(m)\},\\
	B=&\{\pi(1),\pi(2),\ldots,\pi(k)\},\\
	A'=&\{\pi(i+m-k),\pi(i+m-k+1),\ldots,\pi(i+m-1)\}\\
	B'=&\{\pi(i),\pi(i+1),\ldots,\pi(i+k-1)\}.
\end{align*}

Suppose that $i+k-1\leq n$. Since $A\cap B'\neq \emptyset$ and $B\cap B'\neq \emptyset$, we have $n-k+m+1\leq i+k-1$ and $i\leq k$. This implies that $n\leq 3k-m-2$, which contradicts $n\geq 3k-m-1$. 

Now suppose that $i+k-1\geq n+1$. By $n\geq 3k-m-1$, $i+m-k\geq n-2k+m+2\geq k+1$. Also, since $i\leq n-m+1$, we have $i+m-1\leq n$. Thus, $B\cap A'=\emptyset$, which contradicts that $\mathcal{F}$ is intersecting.

Thus, $\mathcal{I}_{m,\pi}$ is an intersecting family consisting of $m$-intervals.
\end{proof}

Then by Lemma~\ref{lem-interval}, we have $|\mathcal{I}_{m,\pi}|\leq m$.

It is clear that there exists a bijection between $\mathcal{I}_{m,\pi}$ and $\mathcal{P}_m\cap \mathcal{R}_{\pi}$. Thus, $|\mathcal{P}_m\cap \mathcal{R}_{\pi}|=|\mathcal{I}_{m,\pi}|\leq m$ and
\begin{align}\label{eq2}
	|X_m|=\sum_{\pi\in C_n}|\mathcal{P}_m\cap \mathcal{R}_{\pi}|\leq \sum_{\pi\in C_n}m= m(n-1)!.
\end{align}

By \eqref{eq2} and \eqref{eq3}, 
\begin{align}\label{eq4}
	P_m(\mathcal{F})\leq \frac{(n-1)!}{(m-1)!(k-m)!^2(n-2k+m)!}=\binom{n-1}{k-1}\binom{k-1}{m-1}\binom{n-k}{k-m}.
\end{align}

By a direct calculation, $P_m(\mathcal{S}_1^k)=\binom{n-1}{k-1}\binom{k-1}{m-1}\binom{n-k}{k-m}$. Thus,
$$P_m(\mathcal{F})\leq P_m(\mathcal{S}_1^k).$$

\qed

\section{Stability}
\subsection{Proof of Theorem~\ref{thm-stability}}

We will give the proof of Theorem~\ref{thm-stability} in this section. If $m=k$, then $n\geq 2k+1$ and Theorem~\ref{thm-stability} follows directly from the Hilton--Milner theorem \cite{Hilton1967}.

Now we assume that $m\in [k-1]$. We still use the notation in Section~3. For fixed $m$ and $\pi$, recall that
$$
\mathcal{I}_{m,\pi}:= \{A\cap B: (A,B)\in \mathcal{P}_m\cap \mathcal{R}_{\pi}\}.
$$
Also, by Claim~\ref{cla-m}, $\mathcal{I}_{m,\pi}$ is an intersecting family consisting of $m$-intervals. By Lemma~\ref{lem-interval}, $|\mathcal{I}_{m,\pi}|\leq m$. We call $\pi$ {\it saturated} if $|\mathcal{I}_{m,\pi}|=m$. Moreover, if $\pi$ is saturated, then all $m$-intervals in $\mathcal{I}_{m,\pi}$ contain a fixed element $x$, which we call the {\it center} of $\mathcal{I}_{m,\pi}$. Then we say $\pi$ is {\it $x$-saturated}. 

We claim that an adjacent transposition does not change the center of saturated cyclic permutations.

\begin{claim}\label{lem-switch}
	Suppose $\pi$ is $x$-saturated and $x=\pi(k)$. For $i\in [n]\backslash \{k-1,k\}$, let $\pi_i\in C_n$ be defined by
	$$
	\pi_i(j)=\left\{
	\begin{aligned}
		\pi(j)\qquad & j\in [n]\backslash\{i,i+1\};\\
		\pi(i+1)\qquad & j=i;\\
		\pi(i)\qquad & j=i+1.
	\end{aligned}
	\right.
	$$
	For $n\geq \max\{2k+1,3k-m-1\}$, if $\pi_i$ is saturated and $(n,m)\neq (3k-2,1)$, then $\pi_i$ is $x$-saturated.
\end{claim}
\begin{proof}[Proof of Claim~\ref{lem-switch}]
    For every $\sigma\in C_n$, let $I_j^{\sigma}$ be the $k$-interval of $\sigma$ with left endpoint $\sigma(j)$. Since $\pi$ is $x$-saturated and $\pi(k)=x$,
    $$
    \mathcal{I}_{m,\pi}=\{\{\pi(k-m+1),\ldots,\pi(k)\},\{\pi(k-m+2),\ldots,\pi(k+1)\},\ldots,\{\pi(k),\ldots,\pi(k+m-1)\}\}.
    $$
    Then
    $$
    \{I_1^{\pi},I_2^{\pi},\ldots,I_m^{\pi}\}\cup\{I_{k-m+1}^{\pi},I_{k-m+2}^{\pi},\ldots,I_{k}^{\pi}\}\subseteq \mathcal{F}.
    $$
    
    We first note that $I_j^{\pi}\notin \mathcal{F}$ if $j\notin [k]$. Indeed, for each $s\in \{1,\ldots,m\}$, $I_s^{\pi}$ and $I_j^{\pi}$ are disjoint when $j\in \{s+k,s+k+1,\ldots,n+s-k\}$. Thus $I_j^{\pi}\notin \mathcal{F}$ for each $j\in \{k+1,k+2,\ldots,n+m-k\}$. Similarly, $I_j^{\pi}\notin \mathcal{F}$ for each $j\in \{2k-m+1,2k-m+2,\ldots,n\}$. Since $n\geq 3k-m-1\geq 3k-2m$, $n+m-k\geq 2k-m$. Thus $I_j^{\pi}\notin \mathcal{F}$ for each $j\in \{k+1,k+2,\ldots,n\}$.
    
    Suppose $\pi_i$ is $\pi_i(c)$-saturated and $c\neq k$. Similarly,
    $$
    \{I_{c-k+1}^{\pi_i},I_{c-k+2}^{\pi_i},\ldots,I_{c-k+m}^{\pi_i}\}\cup\{I_{c-m+1}^{\pi_i},I_{c-m+2}^{\pi_i},\ldots,I_{c}^{\pi_i}\}\subseteq \mathcal{F},
    $$
    where the indices are taken modulo $n$.
    
    {\bf Case 1.} $k+1\leq c\leq 2k-m$.
    
    It is clear that $I^{\pi_i}_{c}= I^{\pi}_{c}$ unless $i=c-1$ or $i=c+k-1$. In this case, since $c\notin [k]$ and $I^{\pi_i}_{c}\in \mathcal{F}$, we have $i=c-1$ or $i=c+k-1$.
    
    If $i=c-1$, then $\pi_i(c)=\pi(c-1)$. Since $i=c-1\neq k$, we have $c\geq k+2$. Moreover, $c+k-1\leq 3k-m-1\leq n$. Thus $I^{\pi_i}_{c}=\{\pi(c-1),\pi(c+1),\ldots,\pi(c+k-1)\}$ is disjoint from  $I^{\pi}_1$, which contradicts that $\mathcal{F}$ is intersecting.
    
    If $i=c+k-1$, then $\pi_i(c+k-1)=\pi(c+k)$. Note that $I^{\pi_i}_{c}\cap I^{\pi}_1=\emptyset$ unless $c+k\equiv 1 \pmod{n}$. Thus $c=n-k+1$ and $i=n$. Since $c\leq 2k-m$ and $n\geq 3k-m-1$, we must have $n=3k-m-1$ and $c=2k-m$. Recall that $(n,m)\neq (3k-2,1)$. So $m\geq 2$ and then $I_{c-1}^{\pi_i}=I_{c-1}^{\pi}\in \mathcal{F}$. However, $n\geq 2k+1$ implies that $c\geq k+2$. Then $c-1\notin [k]$ and $I_{c-1}^{\pi_i}\cap I_1^{\pi}=\emptyset$, which contradicts that $\mathcal{F}$ is intersecting.
    
    {\bf Case 2.} $m\leq c\leq k-1$.
    
    This case is symmetric to Case 1. In this case, $I^{\pi_i}_{n+c-k+1}=I^{\pi_i}_{c-k+1}\in \mathcal{F}$. Similarly, we can prove that $I^{\pi_i}_{n+c-k+1}\cap I^{\pi}_k=\emptyset$ unless $n=3k-m-1$ , $c=m$ and $i=2k-1$. Recall that $(n,m)\neq (3k-2,1)$. So $m\geq 2$ and $I^{\pi_i}_{n+c-k+2}=I^{\pi_i}_{2k+1}=I^{\pi}_{2k+1}\in \mathcal{F}$. However, $n\geq 2k+1$ implies $2k+1\notin [k]$. Hence $I_{2k+1}^{\pi_i}\cap I_k^{\pi}=\emptyset$, which contradicts that $\mathcal{F}$ is intersecting.
    
    {\bf Case 3.} $2k-m+1\leq c\leq n$.
    
    Note that both $I_{c-k+m}^{\pi_i}$ and $I_{c}^{\pi_i}$ are in $\mathcal{F}$. Since $c-k+m\notin [k]$ and $c\notin [k]$, we have $I_{c-k+m}^{\pi_i}\neq I_{c-k+m}^{\pi}$ and $I_{c}^{\pi_i}\neq I_c^{\pi}$. It is clear that $I_{c-k+m}^{\pi_i}\neq I_{c-k+m}^{\pi}$ implies $i\equiv c-k+m-1 \pmod{n}$ or $i\equiv c+m-1 \pmod n$, while  $I_{c}^{\pi_i}\neq I_c^{\pi}$ implies $i\equiv c-1 \pmod{n}$ or $i\equiv c+k-1 \pmod n$. Thus, at least one of
    \begin{align*}
    	c-k+m-1\equiv c-1 \pmod{n},&& c-k+m-1\equiv c+k-1 \pmod{n},\\
    	 c+m-1\equiv c-1 \pmod{n},&& c+m-1\equiv c+k-1 \pmod{n}
    \end{align*}
    holds. Together with $m\in [k-1]$ and $n\geq 2k+1$, this is impossible.
    
    {\bf Case 4.} $1\leq c\leq m-1$.
    
    Note that both $I_{n+c-k+1}^{\pi_i}$ and $I_{n+c-m+1}^{\pi_i}$ are in $\mathcal{F}$. Applying the same argument as in Case~3 gives a contradiction.
    
    Hence, $c=k$, $x=\pi(k)=\pi_i(c)$ and $\pi_i$ is also $x$-saturated.
\end{proof}
    
    We now prove the stability theorem. For a fixed $\varepsilon>0$, choose $\delta<\min\{\frac{1}{m(n^4+1)},\frac{\varepsilon}{m(n^3+1)}\}$. According to \eqref{eq3},
    $$
    \sum_{\pi \in C_n}|\mathcal{I}_{m,\pi}|=\sum_{\pi \in C_n}|\mathcal{P}_m\cap \mathcal{R}_{\pi}|=|X_m|=m!(k-m)!^2(n-2k+m)!P_m(\mathcal{F}).
    $$
    Since $P_m(\mathcal F)\geq (1-\delta)P_m(\mathcal S_1^k)$, 
    \begin{align*}
    	\sum_{\pi \in C_n}|\mathcal{I}_{m,\pi}|&\geq (1-\delta)m!(k-m)!^2(n-2k+m)!P_m(\mathcal S_1^k)\\
    	&=(1-\delta)m(n-1)!,
    \end{align*}
    where the equality holds because $P_m(\mathcal{S}_1^k)=\binom{n-1}{k-1}\binom{k-1}{m-1}\binom{n-k}{k-m}$. Then we have
    $$
    \sum_{\pi \in C_n}\left(m-|\mathcal{I}_{m,\pi}|\right)\leq \delta m(n-1)!.
    $$
    
    Let $U$ be the set of unsaturated cyclic permutations. Since $m-|\mathcal{I}_{m,\pi}|\geq 1$ for each $\pi \in U$, 
    \begin{align}\label{equ}
    	|U|\leq \delta m(n-1)!.
    \end{align}
    
    For each $x\in [n]$, let $S_x$ be the set of all $x$-saturated cyclic permutations. It is clear that the center of every saturated cyclic permutation is unique. Hence $\cup_{x\in [n]}S_x$ is a partition of $C_n\backslash U$.
    
    Now we choose $x\in [n]$ such that $|S_x|$ is maximum. Relabeling the ground set if necessary, assume that $x=k$. For each cyclic order, choose its unique representative $\pi$ satisfying $\pi(k)=k=x$. In this way, $C_n$ is naturally identified with $\mathrm{Sym}([n]\backslash\{k\})$. Let $G_k$ be the Cayley graph defined in Lemma~\ref{lem-cay}. Then each cyclic permutation corresponds to a vertex of $G_k$ and every edge of $G_k$ corresponds to an adjacent transposition which does not involve the position $k$.
    
    By Claim~\ref{lem-switch}, if $\pi\in S_x$ and $\mu$ is adjacent to $\pi$ in $G_k$, then either $\mu\in S_x$ or $\mu \in U$. This implies that $N_{G_k}(S_x)\subseteq U$. 
    
    Recall that $|S_x|$ is maximum. We have 
    $$
    |S_x|\geq \frac{|C_n\backslash U|}{n}=\frac{(n-1)!-|U|}{n}.
    $$
    
    Then we claim that $|S_x|>\frac{|V(G_k)|}{2}$. Suppose otherwise. By Lemma~\ref{lem-cay}, $|N_{G_k}(S_x)|\geq \frac{|S_x|}{n^3}$. Since $N_{G_k}(S_x)\subseteq U$,
    $$
    |U|\geq |N_{G_k}(S_x)|\geq \frac{|S_x|}{n^3}\geq \frac{(n-1)!-|U|}{n^4}.
    $$
    Then $|U|\geq \frac{(n-1)!}{n^4+1}$. Together with \eqref{equ},
    $$
    \delta m(n-1)!\geq |U|\geq \frac{(n-1)!}{n^4+1}.
    $$
    This implies that $\delta\geq \frac{1}{m(n^4+1)}$, which contradicts the choice of $\delta$. Hence $|S_x|>\frac{|V(G_k)|}{2}$.
    
    Let $T$ be the set of all saturated cyclic permutations whose center is different from $x$. Then $|T|\leq \frac{|V(G_k)|}{2}$. By Lemma~\ref{lem-cay}, $|N_{G_k}(T)|\geq \frac{|T|}{n^3}$. Since $N_{G_k}(S_x)\subseteq U$, there is no edge between $S_x$ and $T$ and then $N_{G_k}(T)\subseteq U$. Thus,
    \begin{align}\label{eqt}
    	|U|\geq |N_{G_k}(T)|\geq \frac{|T|}{n^3}.
    \end{align}
    
    By \eqref{equ} and \eqref{eqt}, we have
    \begin{align}
    	|S_x|=&(n-1)!-|U|-|T|\notag\\
    	\geq & (n-1)!-(n^3+1)|U|\notag\\
    	=& (1-(n^3+1)\delta m)(n-1)!.\label{eqs}
    \end{align}
    
    It remains to convert the lower bound on $|S_x|$ into a lower bound on $d_{\mathcal{F}}(x)$. For each $\pi\in S_x$, since $\pi$ is $x$-saturated and $\pi(k)=x$,
    $$
    \mathcal{I}_{m,\pi}=\{\{\pi(k-m+1),\ldots,\pi(k)\},\{\pi(k-m+2),\ldots,\pi(k+1)\},\ldots,\{\pi(k),\ldots,\pi(k+m-1)\}\}.
    $$
    Then
    $$
    \mathcal{I}_{\pi}:=\{I_1,I_2,\ldots,I_m\}\cup\{I_{k-m+1},I_{k-m+2},\ldots,I_{k}\}\subseteq \mathcal{F}_{\pi},
    $$
    where $I_i$ is the $k$-interval of $\pi$ with left endpoint $\pi(i)$. Let $q:=\min\{2m,k\}$. Then $\mathcal{F}_{\pi}$ contains $q$ distinct $k$-intervals $\mathcal{I}_{\pi}$, all of which contain $x$.
    
    We double count the size of the following set 
    $$
    Y_m:=\{(\pi,A):\pi\in S_x,A\in \mathcal{I}_{\pi}\}.
    $$
    
    For a fixed $\pi\in S_x$, there are exactly $q$ intervals belonging to $\mathcal{I}_{\pi}$. Hence, $|Y_m|=q|S_x|.$ On the other hand, for each fixed $A\in \mathcal{F}$ with $x\in A$, there are exactly $(k-1)!(n-k)!$ cyclic permutations such that $A$ is an interval whose left endpoint occurs at a fixed position. Since $A\in \mathcal{I}_{\pi}$, the left endpoint of $A$ occurs at $q$ different positions. Thus,
    $$
    |Y_m|\leq q(k-1)!(n-k)!d_{\mathcal{F}}(x).
    $$
    
    Thus
    \begin{align*}
    	d_{\mathcal{F}}(x)\geq &\frac{|Y_m|}{q(k-1)!(n-k)!}\\
    	=&\frac{|S_x|}{(k-1)!(n-k)!}\\
    	\geq &\frac{(1-(n^3+1)\delta m)(n-1)!}{(k-1)!(n-k)!}\\
    	=&(1-(n^3+1)\delta m)\binom{n-1}{k-1},
    \end{align*}
    where the second inequality holds by \eqref{eqs}. Recall that $\delta<\frac{\varepsilon}{m(n^3+1)}$. We have 
    $d_{\mathcal F}(x)\geq
    (1-\varepsilon)\binom{n-1}{k-1},$ as desired.

\subsection{Proof of the Equality Case in Theorem~\ref{thm-pm}}
    We are now ready to complete the proof of Theorem~\ref{thm-pm}.
	Suppose that $P_m(\mathcal F)=P_m(\mathcal S_1^k)$.
	For any $\varepsilon>0$, Theorem~\ref{thm-stability} gives an $x\in[n]$ such that
	$$
	d_{\mathcal F}(x)\ge(1-\varepsilon)\binom{n-1}{k-1}.
	$$
	Choose $0<\varepsilon<\binom{n-1}{k-1}^{-1}$. Since $d_{\mathcal F}(x)$ is an integer and $d_{\mathcal F}(x)\le\binom{n-1}{k-1}$,
	we obtain
	$$
	d_{\mathcal F}(x)=\binom{n-1}{k-1}.
	$$
	Then $|\mathcal{F}|\geq \binom{n-1}{k-1}$. By Theorem~\ref{thm-ekr}, $\mathcal{F}=\mathcal{S}^k_1$ up to permutations.
    
\section{Concluding remarks}

We conclude with several remarks on the range of the parameters and on some special choices of the weight function $f$.

{\bf \noindent Remark 1.} In Theorem~\ref{thm-pm}, we proved $P_m(\mathcal F)\leq P_m(\mathcal S_1^k)$ when $n\geq\max\{2k,3k-m-1\}$. When $m=k$, Theorem~\ref{thm-pm} reduces to the Erd\H{o}s--Ko--Rado theorem, and hence the condition $n\geq2k$ cannot be weakened. For $m\leq k-1$, the term $3k-m-1$ also cannot, in general, be replaced by $3k-m-2$. The following example shows this.

{\bf Example 1.} Let $m=1, k=3$ and $n=3k-m-2=6\geq 2k$. Let
\begin{align*}
	\mathcal{F}=\{&\{1,2,3\},\{1,2,4\},\{1,3,5\},\{1,4,6\},\{1,5,6\},\\
	&\{2,3,6\},\{2,4,5\},\{2,5,6\},\{3,4,5\},\{3,4,6\}\}.
\end{align*}
Then $P_1(\mathcal{F})=60$. However, $P_1(\mathcal{S}_1^k)=30$. Hence $P_1(\mathcal{F})>P_1(\mathcal{S}_1^k)$.

Moreover, we give a family of counterexamples showing that, for every fixed $m$, the condition $n>2k$ alone is not sufficient for a star to maximize $P_m$.

{\bf Example 2.} For any positive integer $m$, let $k=2m+3$, $n=4m+7$. Let
$$
\mathcal{F}=\left\{A\in \binom{[n]}{k}:|A\cap \{1,2,3\}|\geq 2\right\}.
$$
Then by a direct calculation,  we have
\begin{align*}
	P_m(\mathcal F)
	=&6\binom{n-3}{k-2}\binom{k-2}{m-1}\binom{n-k-1}{k-m-1}+3\binom{n-3}{k-2}\binom{k-2}{m-2}\binom{n-k-1}{k-m}\\
	&+6\binom{n-3}{k-2}\binom{k-2}{m-2}\binom{n-k-1}{k-m-1}+\binom{n-3}{k-3}\binom{k-3}{m-3}\binom{n-k}{k-m}.
\end{align*}
However, $P_m(\mathcal{S}_1^k)=\binom{n-1}{k-1}\binom{k-1}{m-1}\binom{n-k}{k-m}$. Thus,
\begin{align*}
	\frac{P_m(\mathcal F)}{P_m(\mathcal S_1^{k})}&=\frac{6k^2-12km+4m^2+3mn-3n+2}{(n-1)(n-2)}\\
	&=\frac{16m^2+45m+35}{16m^2+44m+30}\\
	&>1.
\end{align*}

{\bf \noindent Remark 2.} In Theorem~\ref{thm-pm}, we also proved that a star is the unique extremal family whenever $n\geq \max\{2k+1,3k-m\}$. Because of the Erd\H{o}s--Ko--Rado theorem, the condition $n\geq 2k+1$ cannot be weakened when $m=k$. For $m\leq k-1$, the condition $n\geq 3k-m$ can be replaced by $n\geq 3k-m-1$ and $(n,m)\neq (3k-2,1)$. This follows because Claim~\ref{lem-switch} holds under $n\geq \max\{2k+1,3k-m-1\}$ with $(n,m)\neq (3k-2,1)$ and the proof of Theorem~\ref{thm-stability} goes through unchanged under the same condition. The following example shows $(n,m)\neq (3k-2,1)$ cannot be omitted for Claim~\ref{lem-switch}.

{\bf Example 3.} Let $m=1, k=4, n=10=3k-2$. Consider the cyclic permutation $\pi=(1,2,\ldots,10)$. Let 
$$
\mathcal{F}=\{\{1,2,3,4\},\{4,5,6,7\},\{1,7,8,9\}\}.
$$
Then $\pi$ is $4$-saturated. However, $\pi_{10}=(10,2,3,\ldots,9,1)$ is $7$-saturated.

The example above shows only that our switching argument fails when $(n,m)=(3k-2,1)$. We conjecture that the uniqueness conclusion of Theorem~\ref{thm-pm} remains true in this boundary case, but a different method may be needed.

{\bf \noindent Remark 3.} Corollary~\ref{cor-phi} shows that a star maximizes $\Phi_f$ for any nonnegative function $f$ when $n\geq 3k-2$. Various choices of the weight function $f$ lead to intersection parameters that have been studied previously. We give some examples in the following.

{\bf Example 4.} Let $f(x)=x$. Then 
$$
\Phi_f(\mathcal{F})=\sum_{(A,B)\in\mathcal F\times\mathcal F}|A\cap B|.
$$
$\Phi_f(\mathcal{F})$ is exactly the {\it total intersection number} of $\mathcal{F}$, which was introduced by Ge and Kong \cite{Ge2026}. It is closely related to the parameter $\omega(\mathcal{F})=\sum_{\{A,B\}\subseteq\mathcal F}|A\cap B|$ studied by Huang, Katona and Yue \cite{Huang2026}. Indeed, $\Phi_f(\mathcal{F})=2\omega(\mathcal{F})+k|\mathcal{F}|$. Note that $\Phi_f(\mathcal{F})=\sum_{x\in [n]}d_{\mathcal{F}}(x)^2$. Hence maximizing $\Phi_f(\mathcal{F})$ is equivalent to maximizing the sum of squares of degrees. Using Bey's result \cite{Bey2003} and the Erd\H{o}s--Ko--Rado theorem, $\sum_{x\in [n]}d_{\mathcal{F}}(x)^2$ is maximized by a star among all intersecting families with $n\geq 2k$. In contrast, Example~2 shows that, for every fixed $m$, there exist $k$ and $n>2k$ for which $P_m(\mathcal F)>P_m(\mathcal S_1^k)$. Thus, when $n>2k$, although a star need not maximize each individual level of the intersection profile, it still maximizes the total intersection number.

{\bf Example 5.} Let $f(x)=\binom{x}{r}$. Then
$$
\Phi_f(\mathcal{F})=\sum_{R\in \binom{[n]}{r}}d_{\mathcal{F}}(R)^2,
$$
where $d_{\mathcal{F}}(R):=|\{A\in \mathcal{F}:R\subseteq A\}|.$ In particular, when $r=k-1$, this is the {\it codegree squared sum} studied recently by Wu and Zhang \cite{Wu2026}.

Although the range of n obtained from our general result is not always as strong as the best known bounds for particular choices of $f$, Corollary~\ref{cor-phi} recovers these earlier results in a common range of parameters.

\section*{Acknowledgements}

The author would like to express his sincere gratitude to Professor Gyula O. H. Katona for his guidance and encouragement. The author had the privilege of spending one year at the Alfr\'ed R\'enyi Institute of Mathematics under Professor Katona's guidance, during which Professor Katona introduced him to the circle method. This method has since played an important role in the author's research and, in particular, in the present work. 

This work was supported by the Startup Foundation for Introducing Talent of NUIST.

\end{document}